\documentclass[10pt]{amsart}

\usepackage[margin=1.5in]{geometry} 
\usepackage{amsmath, amsfonts, amsthm, amstext, amssymb} 
\usepackage{hyperref} 
\usepackage{todonotes} 
\usepackage[nameinlink, capitalize]{cleveref} 
\usepackage{tikz-cd} 
\usepackage[shortlabels]{enumitem} 
\usepackage{array}
\usepackage{scalerel}
\usepackage{mathtools}
\usepackage{cancel} 
\usepackage{lscape} 
\usepackage{fancyhdr} 

\numberwithin{equation}{section} 
\numberwithin{table}{section} 

\theoremstyle{definition}
\newtheorem{definition}[equation]{Definition}

\newtheorem{example}[equation]{Example}

\newtheorem{construction}[equation]{Construction}
\theoremstyle{plain}
\newtheorem{lemma}[equation]{Lemma}
\newtheorem{proposition}[equation]{Proposition}

\newtheorem{maintheorem}{Theorem} 

\def\bbZ{\mathbb{Z}}
\def\z{\mathbb{Z}}

\def\uA{\underline{A}}
\def\uB{\underline{B}}
\def\uC{\underline{C}}
\def\uF{\underline{F}}
\def\ug{\underline{g}}

\def\uK{\underline{K}}

\def\uL{\underline{L}}
\def\uM{\underline{M}}
\def\uN{\underline{N}}

\def\uR{\underline{R}}

\def\uT{\underline{T}}

\def\uTor{\underline{\operatorname{Tor}}}

\DeclareMathOperator{\id}{id} 
\DeclareMathOperator{\Tor}{Tor}

\DeclareMathOperator{\res}{res}
\DeclareMathOperator{\tr}{tr}
\DeclareMathOperator{\nm}{nm}

\makeatletter
\newcommand*{\boxwedge}{\mathbin{\mathpalette\@boxwedge{}}}
\newcommand*{\@boxwedge}[2]{%
  \sbox0{$#1\boxplus\m@th$}%
  \dimen2=.5\dimexpr\wd0-\ht0-\dp0\relax 
  \dimen@=\dimexpr\ht0+\dp0\relax
  \def\lw{.07}
  \kern\dimen2 
  \tikz[
    line width=\lw\dimen@,
    line join=round,
    x=\dimen@,
    y=\dimen@,
  ]
  \draw
    (\lw/2,0) rectangle (1-\lw,1-\lw)
    (\lw,0) -- (.5,1-\lw-\lw/2) -- (1-\lw-\lw/2 ,0)
  ;%
  \kern\dimen2 
}
\makeatother

\newcommand{\ssfrac}[2]{\mathchoice{\raisebox{5pt}{$#1$}\!\big/\!\raisebox{-5pt}{$#2$}}{{}^{#1}\!/_{\!#2}}{ #1 /  #2}{ #1 /  #2}}

\newcommand{\bigsfrac}[2]{\displaystyle \ssfrac{#1}{#2}}

\begin{document}

\author{David Mehrle}\address{Department of Mathematics, University of Kentucky, Lexington, Kentucky, U.S.A.}\email{davidm@uky.edu}
\author{J.D. Quigley}\address{Department of Mathematics, University of Virginia, Charlottesville, VA, U.S.A.}\email{mpb6pj@virginia.edu}
\author{Michael Stahlhauer}\address{Department of Mathematics, Universit\"at Bonn, Bonn, Germany}\email{mstahlhauer@uni-bonn.de}
\title[Pathological computations of Mackey functor-valued Tor]{Pathological computations of Mackey functor-valued Tor over cyclic groups}
\maketitle
\thispagestyle{fancy}
\lfoot{\raisebox{-3em}{2020 Mathematics Subject Classification 18G15, 55P91, 19A22}}
\cfoot{}
\rhead{}
\renewcommand{\headrulewidth}{0pt}

\begin{abstract}
In equivariant algebra, Mackey functors play the role of abelian groups and Green and Tambara functors play the role of commutative rings. In this paper, we compute Mackey functor-valued Tor over certain free Green and Tambara functors, generalizing the computation of Tor over a polynomial ring on one generator. In contrast with the classical situation where the resulting Tor groups vanish above degree one, we present examples where Tor is nonvanishing in almost every degree.
\end{abstract}

\setcounter{tocdepth}{1}
\tableofcontents

\section{Introduction}

Equivariant stable homotopy theory for cyclic $2$-groups has had several remarkable applications in recent years, such as the Hill--Hopkins--Ravenel solution to the Kervaire invariant one problem using $C_8$-equivariant homotopy theory \cite{HHR2016} and the study of the stable homotopy groups of spheres at chromatic height four by Hill--Shi--Wang--Xu using $C_4$-equivariant homotopy theory \cite{HSWX23}. These results involve extensive computations using spectral sequences valued in \emph{Mackey functors} and \emph{Green and Tambara functors}, which are the equivariant analogues of abelian groups and commutative rings, respectively. 

In this paper, we present two computations involving Mackey, Green, and Tambara functors over cyclic groups of prime order. On one hand, our goal is to demonstrate some new techniques (e.g., Koszul-like resolutions) for performing algebra with these highly structured objects. On the other hand, we wish to show just how poorly behaved certain equivariant algebraic constructions can be, similar to the theme of \cite{HMQ21}. 

Before proceeding to our results, we recall (see \cite[Secs. 2-3]{MQS24a} for further details) that over $G = C_p$ a cyclic group of prime order, a Mackey functor $\uM$ (resp. Green functor $\uR$) is a pair of abelian groups  $\uM(C_p/e)$ and $\uM(C_p/C_p)$ (resp. commutative rings $\uR(C_p/e)$ and $\uR(C_p/C_p)$), together with structure maps (restriction, transfer, and conjugation) between them, satisfying certain axioms. A Tambara functor is a Green functor equipped with one additional structure map (the norm) satisfying some additional axioms. There is a well-defined analogue of the tensor product (the box product) which allows for a theory of modules over Green and Tambara functors, and just as in classical algebra, one can define Mackey functor-valued $\uTor$ as the derived functors of the box product. 

In this paper, we focus on computing Mackey functor-valued $\uTor$ over free Green and Tambara functors, which is the equivariant analogue of computing $\Tor$ over polynomial rings. Our motivation for studying $\Tor$ in this equivariant context is its connection to equivariant stable homotopy theory and its role in a conjectural analogue of the Hochschild--Kostant--Rosenberg theorem for Green and Tambara functors. We refer to our companion paper \cite[Sections 1, 4]{MQS24a} for further discussion. 

In ordinary algebra, $\Tor$ over a polynomial ring on one generator is quite simple: 
\[
	\Tor^{\mathbb{Z}[x]}_i(\mathbb{Z},\mathbb{Z}) \cong 
	\begin{cases}
		\mathbb{Z} \quad & \text{ if } i=0,1, \\
		0 \quad & \text{ otherwise.}
	\end{cases}
\]
Our theorem says that Mackey functor-valued $\uTor$ over certain equivariant analogues of polynomial rings on one generator can be infinitely more complicated. In the theorem, $\uA$ denotes the Burnside Mackey functor (the equivariant analogue of $\mathbb{Z}$), $\uA^{\bot}[x_{C_p/e}]$ denotes the free Green functor on an underlying generator (a Green functor analogue of $\mathbb{Z}[x]$), and $\uA^{\top}[x_{C_p/C_p}]$ denotes the free Tambara functor on a fixed generator (a Tambara functor analogue of $\mathbb{Z}[x]$). These objects are described explicitly in \cref{freeGreen} and \cref{lemma:free_Tambara_on_fixed}. 

\begin{maintheorem} \ 
Let $p$ be a prime and let $C_p$ denote the cyclic group of order $p$. 
	\begin{enumerate}[(a)]

		\item Let $\uR = \uA^\bot[x_{C_p/e}]$ be a free $C_p$-Green functor on an underlying generator. The graded Mackey functor $\uTor_*^{\uR}(\uA, \uA)$ is computed in \cref{Prop:CpGreenUnderlying}.  In particular, $\uTor_i^{\uR}(\uA,\uA)$ has infinite rank for almost all $i \geq p$. \\

		\item Let $\uR = \uA^\top[x_{C_p/C_p}]$ be a free $C_p$-Tambara functor on a fixed generator. The graded Mackey functor $\uTor_*^{\uR}(\uA, \uA)$ is computed in \cref{proposition:FreeCpTambaraFixed}. In particular, $\uTor_i^{\uR}(\uA,\uA)$ has finite rank for all $i \geq 0$, but its rank tends to infinity as $i$ tends to infinity.
		
	\end{enumerate}

\end{maintheorem}

This exotic homological behavior is predicted by \cite[Table A.1]{HMQ21}, which shows that the underlying Mackey functors of $\uA^{\bot}[x_{C_p/e}]$ and $\uA^{\top}[x_{C_p/C_p}]$ are not flat. In \cite{MQS24a}, we show that Mackey-valued $\uTor$ is much simpler over free Green and Tambara functors with flat underlying Mackey functors; the reader is invited to compare the computations above with \cite[Prop. 5.1]{MQS24a} and \cite[Thm. 5.3]{MQS24a}. 




\subsection{Acknowledgements}

The authors thank the referee for detailed comments and suggestions. The first and second authors were supported by NSF RTG grant DMS-2135884. 
The second author was supported by an AMS-Simons Travel Grant and NSF grants DMS-2039316 and  DMS-2414922 (formerly DMS-2203785 and DMS-2314082). 
The second and third authors were supported by the Max Planck Institute for Mathematics in Bonn. 


\section{Background}

We begin with a rapid recollection of the concepts from equivariant algebra involved in our computations. A more detailed introduction appears in \cite[Secs. 1-3]{MQS24a}. The statements in this section are specific to Mackey, Green, and Tambara functors over cyclic groups $C_p$.

\subsubsection{Mackey functors}

A \emph{Mackey functor} \cite{Gre1971} $\uM$ is a pair of abelian groups $\uM(C_p/C_p)$ and $\uM(C_p/e)$ together with a restriction homomorphism 
\[
	\res\colon \uM(C_p/C_p) \to \uM(C_p/e),
\] 
a transfer homomorphism 
\[
	\tr\colon \uM(C_p/e) \to \uM(C_p/C_p),
\] 
and a $C_p$-action on $\uM(C_p/e)$. The restriction and transfer must be equivariant for the $C_p$-action (where $\uM(C_p/C_p)$ is given a trivial action) and satisfy a double coset formula: 
\[
	\res \tr(x) = \sum_{g \in C_p} g \cdot x \quad \text{ for all } x \in \uM(C_p/e).
\]  
See \cite{Web2000}, \cite{Maz2013}, or \cite[Def. 2.1]{MQS24a} for more details.
The data of a Mackey functor can be displayed using a Lewis diagram (as introduced in \cite{Lew1988}):
	\[ 
		\uM = \hspace*{1em}
		\begin{tikzcd}[row sep=large] 
			\uM(C_p/C_p) \ar[d, bend right=50,swap,"\res"] \\
			\uM(C_p/e).
				\ar[u, bend right=50,swap,"\tr"]
				\arrow[out=240,in=300,loop,swap,looseness=3,"C_p"]
		\end{tikzcd}		
	\]
	
\subsection{Green and Tambara functors}

A \emph{Green functor} \cite{Dress1971} $\uR$ is a Mackey functor in which $\uR(C_p/C_p)$ and $\uR(C_p/e)$ are commutative rings, restriction is a ring homomorphism, and transfer is an $\uR(C_p/C_p)$-module homomorphism (where $\uR(C_p/e)$ becomes an $\uR(C_p/C_p)$-module via the restriction). See \cite{Lew1988}, \cite{Maz2013}, or \cite[Def. 2.4]{MQS24a} for more details.

\begin{example}
The \emph{Burnside Green functor} $\uA$ is the Green functor with Lewis diagram
\[ 
		\uA = \hspace*{1em}
		\begin{tikzcd}[row sep=large] 
			\mathbb{Z}[t]/(t^2-pt) \ar[d, bend right=50, swap,"t \mapsto p"] \\
			\mathbb{Z}
				\ar[u, bend right=50,swap,"\cdot t"]
				\arrow[out=240,in=300,loop,swap,looseness=3,"\id"]
		\end{tikzcd}		
	\]
	It is called the Burnside Green functor because $\uA(C_p/H)$ is the Burnside ring $A(H)$. The Burnside functor plays the role of $\bbZ$ in equivariant algebra.
\end{example}

A \emph{Tambara functor} \cite{Tam1993} $\uT$ is a Green functor with one additional structure map 
\[
	\nm: \uT(C_p/e) \to \uT(C_p/C_p)
\] 
called the norm, which is a multiplicative (but not necessarily additive) homomorphism satisfying a multiplicative double coset formula 
\[
	\res \nm(x) = \prod_{g \in C_p} g \cdot x \quad \text{ for all } x \in \uT(C_p/e)
\]
and a distributive law for the norm of a sum (see \cite[Cor. 2.6]{HM2019}).  For details, the reader is referred to \cite{Maz2013} or \cite[Def. 2.8]{MQS24a} and the discussion following it. 

\begin{example}
	The Burnside Green functor becomes a Tambara functor when endowed with the norm 
	\[
		\nm(a) = a + \frac{a^p - a}{p} t.
	\]
\end{example}

\subsection{Modules}

The category of Mackey functors is symmetric monoidal with tensor product given by the \emph{box product}, denoted $\boxtimes$ \cite{Lew1981}. If $\uM$ and $\uN$ are Mackey functors, then $\uM \boxtimes \uN$ is a Mackey functor with Lewis diagram
\[ 
		\uM \boxtimes \uN = \hspace*{1em}
		\begin{tikzcd}[row sep=large] 
			\bigg( \uM(C_p/C_p) \otimes \uN(C_p/C_p) \oplus \uM(C_p/e) \oplus \uN(C_p/e) \bigg) / \sim \ar[d, bend right=50,swap,"\res"] \\
			\uM(C_p/e) \otimes \uN(C_p/e),
				\ar[u, bend right=50,swap,"\tr"]
				\arrow[out=240,in=300,loop,swap,looseness=3,"C_p"]
		\end{tikzcd}		
	\]
where we refer to \cite[Def. 3.1]{HM2019} for the definition of $\sim$ and the structure maps. We note that the Burnside Mackey functor $\uA$ is the monoidal unit. 

Green functors are the commutative monoids for this monoidal structure, and a \emph{module} $\uM$ over a Green or Tambara functor $\uR$ is a module over this commutative monoid. Explicitly, an $\uR$-module is a Mackey functor in which $\uM(C_p/C_p)$ is an $\uR(C_p/C_p)$-module, $\uM(C_p/e)$ is an $\uR(C_p/e)$-module, and the structure maps satisfy certain axioms (see \cite[Def. 2.18]{MQS24a} and the discussion following it). The norm of a Tambara functor is ignored for the definition of modules.

The theory of modules over Green and Tambara functors closely parallels that of modules over commutative rings. The category of $\uR$-modules is an abelian category with enough projectives and injectives, so we may do homological algebra there. In particular, we obtain a Mackey functor-valued $\uTor$ as the left-derived functor of the box product (\cite[Def. 2.27]{MQS24a}) which can be computed via projective resolutions (\cite[Def. 2.26]{MQS24a}). 

\subsection{Free modules}
\label{subsection: free modules}

The \emph{free Mackey functor on a fixed generator}, denoted $\uA\{x_G\}$, is simply $\uA$. The \emph{free Mackey functor on an underlying generator}, denoted $\uA\{x_e\}$, is the Mackey functor with Lewis diagram
	\[ 
		\uA\{x_e\} = \hspace*{1em}
		\begin{tikzcd}[row sep=large] 
			\mathbb{Z} \ar[d, bend right=50,swap,"\Delta"] \\
			\mathbb{Z}^{\oplus p}.
				\ar[u, bend right=50,swap,"\nabla"]
				\arrow[out=240,in=300,loop,swap,looseness=3,"C_p"]
		\end{tikzcd}		
	\]
where $\Delta$ is the diagonal, $\nabla$ is the fold map, and $C_p$ acts by cyclic permutation. In this representation, the group $\uA\{x_e\}(C_p/C_e)$ is spanned by $x_e$ and its Weyl conjugates. These modules are free in the sense that they represent evaluating at the level containing the generator. If $\uR$ is a Green or Tambara functor, the \emph{free $\uR$-module on a generator at level $G/H$} is $\uR\{x_H\} := \uR \boxtimes \uA\{x_H\}$. See \cite[Def. 2.19]{MQS24a} for further discussion. If $H=G$ then this identifies with $\uR$, while if $H=e$, we can identify (cf. \cite[Ex. 2.21]{MQS24a})
	\[ 
		\uR\{x_e\} = \hspace*{1em}
		\begin{tikzcd}[row sep=large] 
			\uR(C_p/C_p)\{\tr(x)\} \ar[d, bend right=50,swap,"\Delta"] \\
			\uR(C_p/e)\{x^{(0)}, x^{(1)}, \ldots, x^{(p-1)}\}
				\ar[u, bend right=50,swap,"\nabla"]
				\arrow[out=240,in=300,loop,swap,looseness=3,"C_p"]
		\end{tikzcd}		
	\]
where $C_p$ acts diagonally on the underlying level by the $C_p$-action on $\uR(C_p/e)$ and by cyclic permutation of the $x^{(i)}$'s. Note that $\tr(x) = \tr(x^{(i)})$ for any $i$. The general transfer is then twisted by the Weyl action to ensure that it is equivariant.

When multiple generators are involved, we use the notation $\uR\{a_H, b_K, c_J, \ldots\}$ for a box product $\uR\{a_H\} \boxtimes \uR\{b_K\} \boxtimes \uR\{c_J\} \boxtimes \cdots$. 

\subsection{Free Green and Tambara functors}

The \emph{free Green functor} (resp. \emph{free Tambara functor}) \emph{on a generator at level $G/H$} is the analogue of a polynomial ring in one generator. We refer to \cite[Def. 3.1]{MQS24a} for a precise definition; we will provide the Lewis diagrams for both cases of interest below. 


\section{Examples with infinite homological dimension}\label{Sec:FirstComputations}

To demonstrate how unwieldy Mackey-valued $\uTor$ can be, we compute $\uTor$ over the free Green functor on an underlying generator (\cref{SSS:CpGreenUnder}) and over the free Tambara functor on a fixed generator (\cref{SSS:CpTambaraFixed}). 

\subsection{$\uTor$ over the free Green functor on an underlying generator}\label{SSS:CpGreenUnder}

We begin by computing $\uTor$ over $\uA[x_e] := \uA^\bot[x_{G/e}]$, the free Green functor on an underlying generator. This Green functor was described explicitly for $p=2$ by Blumberg and Hill in \cite[Lemma 3.2]{BH2019}; a mild generalization of their arguments shows:

\begin{proposition}\label{freeGreen}
The free Green functor on an underlying generator $\uA[x_e]$ has Lewis diagram
	\[ 
		\uA[x_e] = \hspace*{1em}
		\begin{tikzcd}[row sep=large] 
			\z[t,t_{\vec{v}}: \vec{v} \in \z_{\geq 0}^{\times p}] / (t^2=pt, t_{\vec{0}} = t, t_{\vec{v}} = t_{\gamma \vec{v}}, t_{\vec{v}} \cdot t_{\vec{w}} = \sum_{g \in C_p} t_{\vec{v} + g \vec{w}}) \ar[d, bend right=50,swap,"\res"] \\
			\z[x^{(i)} : 0 \leq i \leq p-1],
				\ar[u, bend right=50,swap,"\tr"]
				\arrow[out=240,in=300,loop,swap,looseness=3,"C_p"]
		\end{tikzcd}		
	\]
where $\gamma$ is a generator for $C_p$. The structure maps are given by 
$$\res_e^{C_p}(t_{\vec{v}}) = \sum_{g \in C_p} x^{g \vec{v}}, \quad \tr_e^{C_p}(x^{\vec{v}}) = t_{\vec{v}}, \quad \text{ and } \quad \gamma \cdot x^{(i)} = x^{(i+1)},$$
where if $\vec{v} = (v_0,\ldots,v_{p-1})$, we write $x^{\vec{v}} = (x^{(0)})^{v_0}\cdots (x^{(p-1)})^{v_{p-1}}$,  $\gamma \vec{v} = (v_1,\ldots,v_{p-1},v_0)$, and $\gamma x^{(i)} = x^{(i+1)}$ with superscript taken mod $p$. 
\end{proposition}

\begin{definition}\label{Def:SomeMFs}
	Let $\uL^{\vee} := \uA\{x_e\} / (\tr_e^G x_e)$ be the quotient of the free Mackey functor on an underlying generator by the submodule generated by the transfer of that generator. Since every element at the fixed level of $\uA\{x_e\}$ is a transfer, the fixed level of $\uL^\vee$ is zero. On the underlying level, we must quotient by restrictions of transfers, which leaves us with the following Mackey functor:
	\[ 
		\uL^\vee = \hspace*{1em}
		\begin{tikzcd}[row sep=large] 
			0 \ar[d, bend right=50] \\
			\bigsfrac{\bbZ\{x_0, x_1, \ldots, x_{p-1}\}}{\sum_{i=0}^{p-1} x_i}
				\ar[u, bend right=50]
				\arrow[out=240,in=300,loop,swap,looseness=3,"C_p"]
		\end{tikzcd}		
	\]
Let $\gamma$ be a generator for $C_p$. The Weyl action on the underlying level is given by $\gamma \cdot x_i = x_{i+1}$, with indices taken mod $p$ so that $\gamma \cdot x_{p-1} = x_0$.

	Let $\ug$ be the unique Mackey functor with $\ug(G/G) = \z/p$ and $\ug(G/e) = 0$. 
\end{definition}

\begin{proposition}\label{Prop:CpGreenUnderlying}
For $p=2$,
\[
	\uTor^{\uA[x_e]}_*(\uA,\uA) \cong 
		\begin{cases}
			\uA\{y_G\} \quad & \text{ if } *=0, \\
			\uA\{z_e\} \quad & \text{ if } *=1, \\
			\uL^{\vee} \quad & \text{ if } *=2, \\
			\bigoplus_{j\geq 1} \ug & \text{ if } *=3,\\
			0 \quad & \text{ if } *= 4,5,\\
			\bigoplus_{i\geq 1, j\geq 0} \ug & \text{ if } *=6,\\
			\bigoplus_{i\geq 0}( \uTor_{k-4} \oplus \uTor_{k-3}) & \textup{ if } *>6.
		\end{cases}
\]
For $p\geq 3$, 
\[
	\uTor^{\uA[x_e]}_*(\uA,\uA) \cong 
		\begin{cases}
			\uA\{y_G\} \quad & \text{ if } *=0, \\
			\bigoplus_{{p\choose *}/p} \uA\{z_e\} \quad & \text{ if } 1\leq *\leq p-1, \\
			\bigoplus_{j\geq 1} \ug & \text{ if } *=p,\\
			0 \quad & \text{ if } *= p+1, p+2,\\
			\bigoplus_{i\geq 1, j\geq 0} \ug & \text{ if } *=p+3,\\
			\bigoplus_{i\geq 0}( \uTor_{k-4} \oplus \uTor_{k-3}) & \textup{ if } *>p+3.
		\end{cases}
\]
\end{proposition}

\subsubsection{Proof for $p=2$} Let $\uR := \uA[x_e]$. Our first task is to define a free $\uR$-module resolution of $\uA$. Lewis diagrams for the two free $\uR$-modules we will use are given below; their structure maps are as explained in \cref{subsection: free modules} above. Recall that $\gamma$ is a generator for $C_p$. 
\[ 
		\uR\{\omega_G\} = \hspace*{1em}
		\begin{tikzcd}[row sep=large] 
			\z[t,t_{\vec{v}}: \vec{v} \in \z_{\geq 0}^{\times p}] / (t^2=pt, t_{\vec{0}} = t, t_{\vec{v}} = t_{\gamma \vec{v}}, t_{\vec{v}} \cdot t_{\vec{w}} = \sum_{g \in C_p} t_{\vec{v} + g \vec{w}}) \{\omega\} \ar[d, bend right=50,swap,"\res"] \\
			\z[x^{(i)} : 0 \leq i \leq p-1]\{\res(\omega)\},
				\ar[u, bend right=50,swap,"\tr"]
				\arrow[out=240,in=300,loop,swap,looseness=3,"C_p"]
		\end{tikzcd}		
	\]
	\[ 
		\uR\{\omega_e\} = \hspace*{1em}
		\begin{tikzcd}[row sep=large] 
			\mathbb{Z}[x^{(i)}: 0 \leq i \leq p-1]\{\tr(\omega)\} \ar[d, bend right=50,swap,"\res"] \\
			\z[x^{(i)} : 0 \leq i \leq p-1]\{\omega^{(0)}, \omega^{(1)}, \ldots, \omega^{(p-1)}\}.
				\ar[u, bend right=50,swap,"\tr"]
				\arrow[out=240,in=300,loop,swap,looseness=3,"C_p"]
		\end{tikzcd}		
	\]

\begin{construction}[Resolution for $C_2$]\label{Constr:C2GreenUnderRes}
We build a resolution of $\uA$ as an $\uR$-module using as elementary building blocks the free $\uR$-modules
\[\begin{array}{l<{\quad} l<{\quad} l<{\quad} l}
\uB_0= \uR\{y_G\}, & \uB_1= \uR\{z_e\}, & \uB_2= \uR\{w_e\}, & \uB_3= \uR\{a_{j,G}\}, \\
\uB_4= \uR\{b_{j, G}\}, & \widetilde{\uB}_4= \uR\{\delta_{j,e}\}, & \uB_5= \uR\{d_{j,e}\}, & \widetilde{\uB}_5= \uR\{\epsilon_{j,e}\}, \\
\uB_6= \uR\{f_{j,e}\}, & \widetilde{\uB}_6=\uR\{\zeta_{i,j,G}\}.
\end{array}\]
Here, $i$ and $j$ run over all nonnegative integers $\mathbb{Z}_{\geq 0}$; more precisely, when we write $\uR\{\omega_{i,e}\}$ or similar, we mean
$$\uR\{\omega_{i,e}\} := \uR\{\omega_{i,e}\}_{i \in \mathbb{Z}_{\geq 0}} = \uR\{\omega_{0,e}, \omega_{1,e}, \omega_{2,e}, \ldots\} = \uR\{\omega_{0,e}\} \boxtimes \uR\{\omega_{1,e}\} \boxtimes \uR\{\omega_{2,e}\} \boxtimes \cdots$$

The differentials between these $\uR$-modules are given by
\[\begin{array}{l<{\quad} l}
\beta_0\colon \uB_0 \to \uA,\, y \mapsto 1, x \mapsto 0
	&\beta_1\colon \uB_1\to \uB_0,\, z^{(0)} \mapsto x^{(0)} R(y),\\
\beta_2\colon \uB_2\to \uB_1,\, w^{(0)} \mapsto x^{(0)}z^{(1)} - x^{(1)}z^{(0)},
	&\beta_3\colon \uB_3\to \uB_2,\, a_{j}\mapsto (x^{(0)} x^{(1)})^j T(w),\\
\beta_4\colon \uB_4\to \uB_3,\, b_{j}\mapsto (t-2)a_j, 
	& \widetilde{\beta}_4\colon \widetilde{\uB}_4\to \uB_3,\, \delta_j^{(0)} \mapsto x^{(0)}x^{(1)} R(a_j)-R(a_{j+1}),\\
\beta_5\colon \uB_5\to \uB_4,\, d_j^{(0)} \mapsto R(b_j), 
	& \widetilde{\beta}_5\colon \widetilde{\uB}_5\to \widetilde{\uB}_4,\, \epsilon_j^{(0)} \mapsto \delta_j^{(0)}-\delta_j^{(1)},\\
\beta_6\colon \uB_6\to \uB_5,\, f_{j}^{(0)} \mapsto d_j^{(0)}-d_j^{(1)},
	& \widetilde{\beta}_6\colon \widetilde{\uB}_6\to \widetilde{\uB}_5,\, \zeta_{i,j}\mapsto (x^{(0)}x^{(1)})^i T(\epsilon_j),
\end{array}\]
where we write $R(-)$ and $T(-)$ as shorthand for $\res$ and $\tr$. 
From these $\uR$-modules and differentials, we define the beginning of a resolution of $\uA$ by free $\uR$-modules as follows:
\[\begin{array}{l<{\quad} l r}
\uF_k = \uB_k, & \partial_k = \beta_k & \textup{for } 0\leq k \leq 3,\\
\uF_k = \uB_k\oplus \widetilde{\uB}_k, & \partial_k = \beta_k + \widetilde{\beta}_k & \textup{for } 4\leq k\leq 6.
\end{array}\]
At this point, we observe that some type of periodicity occurs: the kernels of $\beta_6$ and $\widetilde{\beta}_6$ are independent of each other and are infinite sums of $\uR$-modules that already occurred as the kernels of the differentials $\beta_2$ and $\beta_3$, respectively. Hence, we may define 
$$\uF_7 =\bigoplus_{i\geq 0}( \uF_3 \oplus \uF_4), \quad \partial_7= \bigoplus_{i\geq 0}(\partial_3+\partial_4).$$ 
In general, we continue by defining 
\[ \uF_k= \bigoplus_{i\geq 0}(\uF_{k-4} \oplus \uF_{k-3})  \quad \textup{and} \quad \partial_k=\bigoplus_{i\geq 0}(\partial_{k-4}+\partial_{k-3}).\]
\end{construction}

It is a straightforward but tedious exercise to verify the following:

\begin{lemma}
The complex $\uF_\bullet$ from \cref{Constr:C2GreenUnderRes} is a free $\uR$-module resolution of $\uA$. 
\end{lemma}

The desired calculation of $\uTor$ for $p=2$ now follows from calculating $H_\ast (\uF_\bullet\boxtimes_{\uR} \uA)$. The differentials of this chain complex are given by
\[\begin{array}{l<{\quad} l<{\quad} l}
\bar{\partial}_0(y)= 0,
	&\bar{\partial}_1(z^{(0)})=0,
	&\bar{\partial}_2(w^{(0)})=0,\\
\bar{\partial}_3(a_0)= T(w),
	&\bar{\partial}_3(a_{j})=0 \textup{ for } j> 0,
	&\bar{\partial}_4(b_{j})=(t-2)a_j,\\ 
\bar{\partial}_4(\delta_j^{(0)})= -R(a_{j+1}),
	&\bar{\partial}_5(d_j^{(0)})= R(b_j), 
	& \bar{\partial}_5(\epsilon_j^{(0)})= \delta_j^{(0)}-\delta_j^{(1)},\\
\bar{\partial}_6(f_{j}^{(0)})=d_j^{(0)}-d_j^{(1)}
	& \bar{\partial}_6(\zeta_{0,j})= T(\epsilon_j),
	& \bar{\partial}_6(\zeta_{i,j})= 0 \textup{ for } i> 0.
\end{array}\]
After this point, the differentials repeat and the homology can be readily computed. This finishes the calculation for $p=2$.

\subsubsection{Proof for $p>2$}

We replace the beginning of the resolution above by a Koszul resolution since we will need to add multiple generators to kill all terms of the form $x^{(i)}z^{(j)}- x^{(j)} z^{(i)}$ in $\uF_1$. On the one hand, these satisfy Koszul-type relations indexed by sets of three elements of the Weyl group, and on the other hand, the transfers of these elements are not sent to 0 anymore. Taking these alterations into account, we define the following Koszul complex:

\begin{construction}[Resolution for $C_p$, $p$ odd]\label{Constr:CpGreenUnderRes}
We begin with $\uK_0=\uB_0 = \uR\{y_G\}$ and $\uK_1= \uB_1=\uR\{z_e\}$, with differentials given by 
\[ \partial_0 \colon \uK_0 \to \uA,\, y \mapsto 1, x \mapsto 0, \quad \partial_1\colon \uK_1\to \uK_0,\, z^{(0)} \mapsto x^{(0)} R(y). \]
Now, we define for $2\leq n\leq p-1$ the $n$-th Koszul module as
\[ \uK_n= \uR\{ (z^{(i_1)}\wedge \ldots \wedge z^{(i_n)})_e \}_{I_n}, \]
a free $\uR$-module on generators indexed by alternating tensors of Weyl-conjugates of $z$. Here, $I_n$ is a set of representatives $\{i_1, \ldots, i_n\}$ of $n$-element subsets of the Weyl group $W_{C_p}(e)\cong \{ 0, \ldots, p-1\}$ under the diagonal action of $W_{C_p}(e)$. 
Note that by the definition of alternating tensors and by the Weyl action, an $\uR(C_p/e)$-basis of $\uK_n(C_p/e)$ is given by alternating tensors $z^{(i_1)}\wedge \ldots \wedge z^{(i_n)}$ for $i_1\leq \ldots \leq i_n$ and \emph{all} $n$-element subsets $\{i_1, \ldots, i_n\} \subset W_{C_p}(e)$. Using this, we define the usual Koszul differential
\[ \partial_n\colon \uK_n\to \uK_{n-1}, \quad z^{(i_1)}\wedge \ldots \wedge z^{(i_n)} \mapsto \sum_{j=1}^n (-1)^{j-1} x^{(i_j)} \cdot z^{(i_1)}\wedge \ldots \wedge \widehat{z^{(i_j)}} \ldots \wedge z^{(i_n)}.\]
Here, we suppress the additional superscript $(0)$ on all generators. 

As it is a classical Koszul complex, $\uK_\bullet \to \uA$ is exact at level $G/e$, and a straightforward calculation proves exactness at level $G/G$. 

Hence, we only need to describe the kernel of the final map $\partial_{p-1}\colon \uK_{p-1}\to \uK_{p-2}$. We know that $\uK_{p-1}$ is generated by the element $z^{(0)}\wedge \ldots \wedge z^{(p-2)}$, and we have
\[ \partial_{p-1}(z^{(0)}\wedge \ldots \wedge z^{(p-2)}) = \sum_{j=0}^{p-2} (-1)^{j} x^{(j)} \cdot z^{(0)}\wedge \ldots \wedge \widehat{z^{(j)}} \wedge \ldots \wedge z^{(p-2)}. \]
As in the classical Koszul complex, the kernel at level $G/e$ is generated by the element
\begin{equation}\label{eq:KoszulDifferentialAsRestriction}
	\sum_{j=0}^{p-1} (-1)^j x^{(j)} \cdot z^{(0)}\wedge \ldots \wedge \widehat{z^{(j)}} \ldots \wedge z^{(p-1)}.
\end{equation}
This agrees with $\res( x^{(p-1)} \tr(z^{(0)}\wedge \ldots \wedge z^{(p-2)}))$, since the sign is exactly the sign of the permutations of the indices induced by the Weyl group action. In fact, we have $\partial_{p-1}(x^{(p-1)} \tr(z^{(0)}\wedge \ldots \wedge z^{(p-2)}))=0$. Moreover, the kernel of $\partial_{p-1}$ at level $G/G$ is of the form
\[ \{ F\cdot x^{(p-1)} \tr(z^{(0)}\wedge \ldots \wedge z^{(p-2)}) \mid F\in\uR(G/e)^{W_G(e)} \}. \]
Direct investigation of Weyl-fixed elements implies that the kernel is generated by the elements
\[ \nm(x)^j\cdot x^{(p-1)} \tr(z^{(0)}\wedge \ldots \wedge z^{(p-2)}) \quad \textup{for } j\geq 0,\]
where $\nm(x)= \prod_{k=0}^{p-1} x^{(k)}$. Hence, we may continue the resolution by
\[ \uB_3\to \uK_{p-1}, \quad a_j\mapsto \nm(x)^j\cdot x^{(p-1)} \tr(z^{(0)}\wedge \ldots \wedge z^{(p-2)}),\]
where $\uB_3$ is as in \cref{Constr:C2GreenUnderRes}. After this point, the resolution continues as in the case $p=2$, replacing all occurrences of $2$ by $p$.
\end{construction}

The preceding discussion shows: 

\begin{lemma}
The complex $\uK_\bullet$ of \cref{Constr:CpGreenUnderRes} is a free $\uR$-module resolution of $\uA$. 
\end{lemma}

To calculate $\uTor$ in this case, we observe that all differentials in the Koszul complex $\uK_\bullet$ contain coefficients of $x^{(i)}$, and hence reduce to 0 after applying $\_\boxtimes_{\uR} \uA$. In contrast with the case $p=2$, the differential $\uB_3\boxtimes_{\uR} \uA\to \uK_{p-1}\boxtimes_{\uR} \uA$ also evaluates to 0. After this point, all differentials behave as in the case $p=2$. This finishes the calculation for $p$ odd. 

\subsection{$\uTor$ over the free Tambara functor on a fixed generator}\label{SSS:CpTambaraFixed}

Let $p$ be any prime and let $G = C_p$. We now consider $\uTor$ over $\uA[x_G] := \uA^\top[x_{G/G}]$, the free Tambara functor on a fixed generator. This Tambara functor was described explicitly for $p=2$ by Blumberg and Hill \cite[Lemma 3.6]{BH2019}, and as in the Green functor case, a mild generalization of their description yields:

\begin{proposition}\label{lemma:free_Tambara_on_fixed}
The free Tambara functor on a fixed generator $\uA[x_G]$ has Lewis diagram
	\[ 
		\uA[x_G] = \hspace*{1em}
		\begin{tikzcd}[row sep=large] 
			\z[t,n,x]/(t^2=pt,tx^p=tn) \ar[d, bend right=50,swap,"\res"] \\
			\z[x].
				\ar[u, bend right=50,swap,"\tr"] \ar[u,swap,"\nm" description]
				\arrow[out=240,in=300,loop,swap,looseness=3,"C_p"]
		\end{tikzcd}		
	\]
The structure maps are given by
$$\res(t) = p, \ \res(x) = x, \ \res(n)=x^p, \quad \tr(x)=tx, \quad \nm(x)=n, \quad \gamma x = x,$$
where $\gamma$ is a generator for $C_p$.
\end{proposition}

%
%
We introduce the following Mackey functor to describe $\uTor$ below:
\begin{definition}
Let $\uL := \uA\{x_G\} / (\res_e^G x_G)$.
\end{definition}

The goal of this subsection is to make the following computation:

\begin{proposition}\label{proposition:FreeCpTambaraFixed}
We have 
\[
	 \uTor^{\uA[x_{G}]}_i(\uA,\uA) \cong 
		\begin{cases}
			\uA \quad & \text{ if } i=0, \\
			\uA \oplus \uL \quad & \text{ if } i=1, \\
			\uL \quad & \text{ if } i=2, \\
			 0 \quad & \text{ if } i=3,4, \\
			\ug \quad & \text{ if } i=5,6, \\
			0 \quad & \text{ if } i=7,8, \\
			\uT_{i-4} \oplus \uT_{i-4} \oplus \uT_{i-5} \quad & \text{ if } i \geq 9.
		\end{cases}
\]
\end{proposition}

Let $\uR:=\uA[x_G]$. As with our computation over $\uA[x_{e}]$, our first task is to define a free $\uR$-module resolution of $\uA$. Lewis diagrams for the two free $\uR$-modules we will need are listed below; the structure maps are as described in the section on free modules above:
\[ 
		\uR\{\omega_G\} = \hspace*{1em}
		\begin{tikzcd}[row sep=large] 
			\z[t,n,x]/(t^2=pt,tx^p=tn)\{\omega\} \ar[d, bend right=50,swap,"\res"] \\
			\z[x]\{R(\omega)\},
				\ar[u, bend right=50,swap,"\tr"] \ar[u,swap,"\nm" description]
				\arrow[out=240,in=300,loop,swap,looseness=3,"C_p"]
		\end{tikzcd}			
	\]
	\[ 
		\uR\{\omega_e\} = \hspace*{1em}
		\begin{tikzcd}[row sep=large] 
			\z[x]\{T(\omega)\} \ar[d, bend right=50,swap,"\res"] \\
			\z[x]\{\omega^{(0)}, \omega^{(1)}, \ldots, \omega^{(p-1)}\}.
				\ar[u, bend right=50,swap,"\tr"]
				\arrow[out=240,in=300,loop,swap,looseness=3,"C_p"]
		\end{tikzcd}		
	\]

\begin{construction}

Let

\[\begin{array}{l<{\quad} l l}
\uF_0 = \uR\{y_G\}, \quad & \uF_1 = \uR\{z_G,w_G\}, \quad & \uF_2 = \uR\{a_G,b_e\}, \\
\uF_3 = \uR\{c_e,d_e\}, \quad & \uF_4 = \uR\{f_e,h_G\}, \quad & \uF_5 = \uR\{m_G, s_G, \delta_G\}, \\
 \uF_6 = \uR\{u_G, v_G, \xi_e, \zeta_e, \theta_G\}, \quad & \uF_7 = \uR\{\alpha_e, \beta_e, \bar{\alpha}_e, \bar{\beta}_e, \epsilon_G, \omega_e\}.
\end{array}\]
The differentials $\phi_k: \uF_k \to \uF_{k-1}$ (we take $\uF_{-1}=\uA$) between these $\uR$-modules are given by
\[\begin{array}{l<{\quad}l l}
\phi_0(y) = 1, \quad & \quad &  \\
 \phi_1(z) = xy, \quad & \phi_1(w) = (x^p-n)y, \\
\phi_2(a) = (x^p-n)z-xw, \quad & \phi_2(b^{(0)})=R(w), \\
 \phi_3(c^{(0)}) = xb^{(0)}-R(a), \quad & \phi_3(d^{(0)}) = b^{(0)}- b^{(1)}, \\
 \phi_4(h) = T(d), \quad & \phi_4(f^{(0)}) = c^{(0)} - c^{(1)} - x d^{(0)}, \\
 \phi_5(m) = x h + T(f), \quad & \phi_5(s) = (x^p-n)h, \quad & \phi_5(\delta) = (t-p)h, \\
 \phi_6(u) = xs - (x^p-n)m, \quad & \phi_6(v) = x \delta - (t-p)m, \quad &  \phi_6(\xi^{(0)}) = R(s),\\
 \phi_6(\zeta^{(0)}) = R(\delta), \quad & \phi_6(\theta) = (t-p)s-(x^p-n)\delta, \\
 \phi_7(\alpha^{(0)}) = R(u) - x\xi^{(0)}, \quad & \phi_7(\bar{\alpha}^{(0)}) = R(v) - x \zeta^{(0)}, \quad  & \phi_7(\beta^{(0)}) = \xi^{(0)}-\xi^{(1)}, \\
 \phi_7(\bar{\beta}^{(0)}) = \zeta^{(0)} -  \zeta^{(1)}, \quad & \phi_7(\epsilon) = (t-p)u - (x^p-n)v - x \theta, \quad & \phi_7(\omega^{(0)}) = R(\theta).
\end{array}\]







Note that there is some repetition once we get to $\uF_6$: the images of $u_G$ and $\xi_e$ (resp. $v_G$ and $\zeta_e$) should be compared with the images of $a_G$ and $b_e$ in the map $\uF_2 \to \uF_1$. Once we get to $\uF_7$, we see that the kernel of this map is isomorphic to a sum of kernels of previous maps: the kernels of the restriction to the $\alpha_e$ and $\beta_e$ summands (resp. $\bar{\alpha}_e$ and $\bar{\beta}_e$ summands) should be compared with the kernel of the map $\uF_3 \to \uF_2$, and similarly for $\epsilon_G$ and $\omega_e$ and the kernel of $\uF_2 \to \uF_1$. 

For $k \geq 8$, let
$$\uF_k = \uF_{k-4} \oplus \uF_{k-4} \oplus \uF_{k-5}$$
with maps $\uF_k \to \uF_{k-1}$ defined as the sum of maps $\uF_{k-4} \to \uF_{k-5}$ and $\uF_{k-5} \to \uF_{k-6}$.


\end{construction}

As above, a direct analysis of kernels and images shows that $\uF_\bullet$ is indeed a resolution. 

The desired calculation of $\uTor$ now follows from calculating $H_*(\uF_\bullet \boxtimes_{\uR} \uA)$; the differentials of this chain complex are given by
\[\begin{array}{l<{\quad} l l}
 \bar{\phi}_0(y) = 0, \\
 \bar{\phi}_1(z) = 0, \quad & \bar{\phi}_1(w)=0, \\
 \bar{\phi}_2(a) = 0, \quad & \bar{\phi}_2(b^{(0)}) = R(w), \\
 \bar{\phi}_3(c) = -R(a), \quad & \bar{\phi}_3(d^{(0)}) = b^{(0)}- b^{(1)}, \\
 \bar{\phi}_4(h) = T(d), \quad & \bar{\phi}_4(f^{(0)}) = c^{(0)} - c^{(1)}, \\
 \bar{\phi}_5(m) = T(f), \quad & \bar{\phi}_5(s) = 0, \quad & \bar{\phi}_5(\delta) = (t-p)h, \\
 \bar{\phi}_6(u) = 0, \quad & \bar{\phi}_6(\xi^{(0)}) = R(s), \quad & \bar{\phi}_6(v) = -(t-p)m, \\
 \bar{\phi}_6(\zeta^{(0)}) = R(\delta), \quad & \bar{\phi}_6(\theta) = (t-p)s, \\
 \bar{\phi}_7(\alpha^{(0)}) = R(u), \quad & \bar{\phi}_7(\bar{\alpha}^{(0)}) = R(v), \quad & \bar{\phi}_7(\beta^{(0)}) = \xi^{(0)}-\xi^{(1)}, \\
 \bar{\phi}_7(\bar{\beta}^{(0)}) = \zeta^{(0)} - \zeta^{(1)}, \quad & \bar{\phi}_7(\epsilon) = (t-p)u, \quad & \bar{\phi}_7(\omega) = R(\theta),
\end{array}\]
and for $k \geq 8$, we have 
$$\uC_k \cong \uC_{k-4} \oplus \uC_{k-4} \oplus \uC_{k-5}$$
with $\uC_k \to \uC_{k-1}$ given by the sum of the maps $\uC_{k-4} \to \uC_{k-5}$ and $\uC_{k-5} \to \uC_{k-6}$. 

\bibliographystyle{alpha}
\bibliography{references}

\end{document}